\documentclass[reqno]{amsart}
\usepackage{hyperref}

\begin{document}
\title[\hfilneg \hfil Unicity of Entire  Functions Concerning their $q-$ Derivatives-Difference-Polynomials]
{Unicity of Entire  Functions Concerning their $q-$ Derivatives-Difference-Polynomials}

\author[XiaoHuang Huang \hfil \hfilneg]
{XiaoHuang Huang}

\address{XiaoHuang Huang: Corresponding author\newline
Department of Mathematics, Shenzhen University, Shenzhen 518055, China}
\email{1838394005@qq.com}


\begin{abstract}
In this paper, we study the unicity of entire functions concerning their $q-$shifts and  $k-$th derivatives and prove:
Let $f(z)$ be a transcendental entire function of zero-order, and $g(z)$ define as in (1.1). Let $a(z), b(z)$ be two distinct  small functions of $f(z)$. If $f(z)$ and $g(z)$ share $a(z), b(z)$ IM, then $f(z)\equiv g(z)$.
\end{abstract}

\maketitle
\numberwithin{equation}{section}
\newtheorem{theorem}{Theorem}[section]
\newtheorem{lemma}[theorem]{Lemma}
\newtheorem{remark}[theorem]{Remark}
\newtheorem{corollary}[theorem]{Corollary}
\newtheorem{example}[theorem]{Example}
\newtheorem{problem}[theorem]{Problem}
\allowdisplaybreaks

\section{Introduction and main results}

Let $\Bbb C$ denote the complex plane and $f$ a meromorphic function on $\Bbb C$. In this paper, we assume that the reader is familiar with the fundamental results
and the standard notation of the Nevanlinna value distribution theory, see(\cite{h3,y1,y2}). In addition, $S(r, f) = o(T(r, f))$, as $r\to \infty $ outside of a possible exceptional set of finite logarithmic measure. Define
 $$\rho(f)=\varlimsup_{r\rightarrow\infty}\frac{log^{+}T(r,f)}{logr},$$
 $$ \mu(f)=\varliminf_{r\rightarrow\infty}\frac{log^{+}T(r,f)}{logr},$$
as the order and the lower order  of $f$.

For a meromorphic function $f(z)$, we define its $q-$shift by $f(qz+c)$.

Let $f$ and $g$ be two nonconstant meromorphic functions, and let $a(z)$ be a small function of $f$ and $g$. We say that $f$ and $g$ share $a$ CM(IM), provided that $f(z)-a$ and $g(z)-a$ have the same zeros counting multiplicities(ignoring multiplicities). Moreover, we introduce the following notation: $S_{(m,n)}(a)=\{z|z $ is a common zero of $f(z)-a$ and $g(z)-a$ with multiplicities $m$ and $n$ respectively$\}$. $\overline{N}_{(m,n)}(r,\frac{1}{f-a})$ denotes the counting function of $f$ with respect to the set $S_{(m,n)}(a)$. $\overline{N}_{n)}(r,\frac{1}{f-a})$ denotes the counting function of all zeros of $f-a$ with multiplicities at most $n$. $\overline{N}_{(n}(r,\frac{1}{f-a})$ denotes the counting function of all zeros of $f-a$ with multiplicities at least $n$. $\overline{N}_{n}(r,\frac{1}{f-a})$ denotes the counting function of all zeros of $f-a$ with multiplicity $n$.

Let $f(z)$ be a meromorphic function, and two finite complex number $q\neq0, c\neq0$, we define  its difference operators by
\begin{equation*}
\Delta_{q,c} f(z)=f(qz+c)-f(z), \quad \Delta^{n}f(z)=\Delta^{n-1}_{q,c}(\Delta_{q,c} f(z)).
\end{equation*}

Rubel and Yang \cite{ruy} first investigated the uniqueness of an entire function concerning its derivative, and proved the following result.

\

{\bf Theorem A} \ Let $f$ be a non-constant entire function, and let $a, b$ be two finite distinct complex values. If $f(z)$ and $f'(z)$
 share $a, b$ CM, then $f(z)\equiv f'(z)$.

Mues and Steinmetz \cite {muess} improved Theorem A and proved

\

{\bf Theorem B} \ Let $f$ be a non-constant entire function, and let $a, b$ be two finite distinct complex values. If $f(z)$ and $f'(z)$
 share $a, b$ IM, then $f(z)\equiv f'(z)$.

In recent years, there has been many interests in developing  the value distribution of meromorphic functions with respect to difference analogue, see [3,4,5,7,8,9,15].
 Heittokangas et al \cite{hkl} proved a similar result analogue of Theorem A concerning shifts. They obtained if a non-constant finite order entire function share two distinct finite values CM with its shift, then they must be identity equal.

Qi \cite {q} improved Theorem B and get a version of Theorem B concerning shifts.

\

{\bf Theorem C}
 Let $f(z)$ be a non-constant entire function of finite order, let $\eta$ be a nonzero finite complex value, and let $a, b$ be two finite distinct complex values.
If $f(z)$ and $f(z+\eta)$ share $a, b$ IM, then $f(z)\equiv f(z+\eta).$

Over last decade, a number of papers concerning $q$-difference and  $q$-shift of a meromorphic function were published, see [2,10,11,13,14]. We are concerned in this paper with respect to $q$-difference and  $q$-shift.  In 2011, Qi-Liu-Yang \cite{qly} proved an analogue of Theorem C.

\

{\bf Theorem D}
 Let $f(z)$ be a non-constant entire function of zero-order, let $q$ be a nonzero finite complex value, and let $a, b$ be two finite distinct complex values.
If $f(z)$ and $f(qz)$ share $a, b$ IM, then $f(z)\equiv f(qz).$

In the following, we define  $g(z)$  to be the $q$-shift differential polynomial of $f(z)$ as
\begin{align}
g(z):=A_{0}(z)f^{(k_{0})}(q_{0}z+c_{0})+A_{1}(z)f^{(k_{1})}(q_{1}z+c_{1})+\cdots+A_{j}(z)f^{(k_{j})}(q_{j}z+c_{j}),
\end{align}
where for $i=0,1,\ldots,j$, $A_{i}(z)$ are some distinct small functions of $f(z)$, $k_{i}$ are distinct positive integers, $c_{i}$ are some distinct finite complex numbers and $q_{i}\neq0$ are some distinct finite complex numbers.

Next, we analyze the order of $g(z)$. By the Lemma 2.1 in the following, we can obtain that for all $i=0,1,\ldots,j$, $\rho(f(z))=\rho(f(q_{i}z+c_{i}))$ and $\rho(f(z))=\rho(f^{k_{i}}(z))$. Since $T(r,A_{i}(z))=o(T(r,f))$, we have $\rho(A_{i}(z))\leq\rho(f(z))$. Therefore, by Theorem 1.16 in \cite{y1}, we have $\rho(g(z))\leq\{\rho(A_{i}(z)),\rho(f^{k_{i}}(q_{i}z+c_{i}))\}=\rho(f(z))$ for all $i=0,1,\ldots,j$. That is to say, if $f(z)$ is a meromorphic function of zero-order, $g(z)$ is also a meromorphic function of zero-order.

It is natural to ask a question that:
\

{\bf Question 1} As in Theorem D, can we replace two distinct finite values by two distinct small function, and replace $f(qz)$ by $g(z)$?

In this paper, we give a positive answer to question 1.  We obtain.

\

{\bf Theorem 1} Let $f(z)$ be a transcendental entire function of zero-order, and $g(z)$ define as in (1.1). Let $a(z), b(z)$ be two distinct  small functions of $f(z)$. If $f(z)$ and $g(z)$ share $a(z), b(z)$ IM, then $f(z)\equiv g(z)$.

Immediately, we have

\

{\bf Corollary} Let $f(z)$ be a transcendental entire function of zero-order, let $n$ be a positive integer, $q\neq0$ a finite complex number, and let  $a(z), b(z)$ be two distinct  small functions of $f(z)$. If $f(z)$ and $\Delta^{n}f(qz)$ share $a(z), b(z)$ IM, then $f(z)\equiv \Delta^{n}_{q,c}f(z)$.

\section{Some Lemmas}
\begin{lemma}\label{21l}\cite{bh} Let $f$ be a nonconstant meromorphic function zero-order,  and let $q$ be a non-zero complex number. Then
$$m(r,\frac{f(qz)}{f(z)})=o(T(r,f)),$$
for all $r$ on a set of logarithmic density $1$.
\end{lemma}

\begin{lemma}\label{23l} Let $f_{1}$ and $f_{2}$ be  nonconstant meromorphic functions in $|z|<\infty$, then
$$N(r,f_{1}f_{2})-N(r,\frac{1}{f_{1}f_{2}})=N(r,f_{1})+N(r,f_{2})-N(r,\frac{1}{f_{1}})-N(r,\frac{1}{f_{2}}),$$
where $0<r<\infty$.
\end{lemma}

\begin{lemma}\label{21l}\cite{bh} Let $f$ be a nonconstant meromorphic function of zero-order,  and let $c$ be a non-zero complex number. Then
$$T(r,f(z))=T(r,f(qz))+o(T(r,f)).$$
\end{lemma}

\begin{lemma}\label{211} \cite{y1} Let $f$ is a non-constant meromorphic function, and let $a_{1}, a_{2}, a_{3}$ be three distinct small functions of $f$. Then
$$T(r,f)\leq \sum_{i=1}^{3}\overline{N}(r,\frac{1}{f-a_{i}})+o(T(r,f)).$$
\end{lemma}

\begin{lemma}\label{23l} Let $f(z)$ be a transcendental entire function, let $k$ be a  positive integer, and let $a(z)\not\equiv\infty, b(z)\not\equiv\infty$ be two distinct small meromorphic functions of $f(z)$. Suppose
\[L(f(z))=\left|\begin{array}{rrrr}a-b& &f(z)-a \\
a'-b'& &f'(z)-a'\end{array}\right|\]
and
\[L(g(z))=\left|\begin{array}{rrrr}a-b& &g(z)-a \\
a'-b'& &g'(z)-a'\end{array}\right|,\]
and $f(z)$ and $g(z)$ share $a$  and  $b$ IM,  then $L(f(z))\not\equiv0$ and $L(g(z))\not\equiv0$.
\end{lemma}
\begin{proof}
Suppose that $L(f(z))\equiv0$, then we can get $\frac{f'(z)-a'}{f(z)-a}\equiv\frac{a'-b'}{a-b}$. Integrating both side of above we can obtain $f(z)-a=C_{1}(a-b)$, where $C_{1}$ is a nonzero constant. Then we have $T(r,f(z))=T(r,f(z))+o(T(r,f))=o(T(r,f))$, a contradiction. Hence $L(f(z))\not\equiv0$.

Since $g(z)$ and $f(z)$ share $a$ and $b$ IM, and that $f(z)$ is a non-constant entire function, then by Lemma 2.1, we get
\begin{align}
T(r,f(z))&\leq \overline{N}(r,\frac{1}{f(z)-a})+\overline {N}(r,\frac{1}{f(z)-b})+o(T(r,f))\notag\\
&= \overline {N}(r,\frac{1}{g(z)-a})+\overline {N}(r,\frac{1}{g(z)-b})+o(T(r,f))\notag\\
&\leq 2T(r,g(z))+o(T(r,f)).
\end{align}
Hence $a$ and $b$ are small functions of $g(z)$. If $L(g(z))\equiv0$, then we can get $g(z)-a=C_{2}(a-b)$, where $C_{2}$ is a nonzero constant. And we get $T(r,g(z))=o(T(r,f))$. Combing (2.1) we obtain $T(r,f(z))=o(T(r,f))$, a contradiction.
\end{proof}

\begin{lemma}\label{24l}  Let $f(z)$ be a transcendental entire function, and $k$ a positive integer. Let $a(z)\not\equiv\infty, b(z)\not\equiv\infty$ be two distinct small meromorphic functions of $f(z)$.  Again let $d_{j}=a-j(a-b)$, where $j\neq0,1$ is a positive integer. Then
$$m(r,\frac{L(f(z))}{f(z)-a})=o(T(r,f)), \quad m(r,\frac{L(f(z))}{f(z)-b})=o(T(r,f)).$$
And
$$m(r,\frac{L(f(z))f(z)}{(f(z)-a)(f(z)-b)(f(z)-d_{j})})=o(T(r,f)),$$
where $L(f(z))$ is defined as in Lemma 2.5.
\end{lemma}
\begin{proof}
Obviously, we have
$$m(r,\frac{L(f(z))}{f(z)-a})\leq m(r,-\frac{(a'-b')(f(z)-a)}{f(z)-a})+m(r,\frac{(a-b)(f'(z)-a')}{f(z)-a})=o(T(r,f)),$$
and
$$\frac{L(f(z))f(z)}{(f(z)-a)(f(z)-b)(f(z)-d_{j})}=\frac{C_{1}L(f(z))}{f(z)-a}+\frac{C_{2}L(f(z))}{f(z)-b}+\frac{C_{3}L(f(z))}{f(qz)-d_{j}},$$
where $C_{i}(i=1,2,3)$ are small functions of $f$. Thus  we have
\begin{align}
m(r,\frac{L(f(z))f(z)}{(f(z)-a)(f(z)-b)(f(z)-d_{j})})=o(T(r,f)).
\end{align}
\end{proof}

\begin{lemma}\label{251}\cite{h3,y1,y2} Suppose that $f(z)$ is a meromorphic function and $p(f)= a_{0}f^{n}(z)+a_{1}f^{n-1}(z)+\cdots +a_{n}$, where $a_{0}(\not\equiv0)$, $a_{1}$,$\cdots$,$a_{n}$ are small functions of $f(z)$. Then
$$T(r,p(f))= nT(r,f(z))+o(T(r,f)).$$
\end{lemma}

In 2013, K. Yamanoi\cite{ya} proved the famous Gol’dberg conjecture, and in his paper, he obtained a more general result.

\begin{lemma}\label{211} \cite{ya} Let $f$ be a transcendental meromorphic function in the complex plane. Let
$k\geq 2$ be an integer, and let $\varepsilon \geq\varepsilon_{1}> 0$. Let $A $ be a finite set of finite complex numbers. Then we have
$$(k-1)\overline{N}(r,f(z))+\sum_{a\in A}N_{1}(r,\frac{1}{f(z)-a})= N(r,\frac{1}{f^{(k)}(z)})+\varepsilon_{1}T(r,f),$$
for all $r>e$ outside a set $E\subset(e,\infty)$  of logarithmic density $0$. Here $E$ depends on $f,A,k$ and $\varepsilon$, and where
$$N_{1}(r,\frac{1}{f(z)-a})=N(r,\frac{1}{f(z)-a})-\overline{N}(r,\frac{1}{f(z)-a}).$$
\end{lemma}

\

{\bf Remark 1} In Lemma 2.7, we set
$$S(r)=(k-1)\overline{N}(r,f(z))+\sum_{a(z)\in A}N_{1}(r,\frac{1}{f(z)-a(z)})-N(r,\frac{1}{f^{(k)}(z)}).$$
With a similar method of proving Lemma 1.5 in \cite{ya}, one can verify that the logarithmic density
of the "exceptional set"
$$E_{\varepsilon}=\{r>e;|S(r)|>\varepsilon T(r,f)\}$$
is zero. That is to say, $|S(r)|\leq\varepsilon T(r,f)$ holds for all $r>e$  outside some exceptional set of logarithmic density zero. And then we can find an $-\varepsilon\leq\varepsilon_{1}\leq\varepsilon$ such that $S(r)=\varepsilon_{1}T(r)$.

Let $a$ be a value in the extended complex plane. We define the deficiency of $a$ with respect to $f(z)$ as
$$\delta(a,f(z))=1-\varlimsup_{r\rightarrow\infty}\frac{N(r,\frac{f-a}{})}{T(r,f)},$$
and if $\delta(a,f(z))>0$, we say that $a$ is a deficient value of $f(z)$.

\begin{lemma}\label{211}\cite{ef}
Meromorphic functions with more than one deficient value have a positive lower order.
\end{lemma}

{\bf Remark 2}\quad If $f(z)$ is a meromorphic function of zero-order with $\delta(\infty,f(z))=1$, then for any finite value $a$, $\delta(a,f(z))=0$ holds.

\begin{lemma}\label{211}\cite{t}
Let $f$ be a transcendental entire function of lower order zero and let $k$ be a positive integer. Then
$$T(r,f(z))\leq T(Kr,f(z))+o(T(r,f))\leq T(r,f^{(k)}(z))+o(T(r,f))\leq T(r,f(z))+o(T(r,f)),$$
for any $K\geq1$.
\end{lemma}

\section{The proof of Theorem 1}
Assume that $f(z)\not\equiv g(z)$. Since $f(z)$ and $g(z)$ share $a$ and $b$ IM, and $f$ is a transcendental entire function of zero-order, then by the Nevanlinna Second Fundamental Theorem, and Lemma 2.1, we get
\begin{eqnarray*}
\begin{aligned}
T(r,f(z))&\leq \overline{N}(r,\frac{1}{f(z)-a})+\overline{N}(r,\frac{1}{f(z)-b})+o(T(r,f))\\
&= \overline{N}(r,\frac{1}{g(z)-a})+\overline{N}(r,\frac{1}{g(z)-b})+o(T(r,f))\\
&\leq N(r,\frac{1}{f(z)-g(z)})+o(T(r,f))\\
&\leq T(r,f(z)-g(z))+o(T(r,f))\\
&\leq m(r,f(z)-g(z))+o(T(r,f))\\
&\leq m(r,f(z))+m(r,1-\frac{g(z)}{f(z)})+o(T(r,f))\\
&\leq T(r,f(z))+o(T(r,f)).
\end{aligned}
\end{eqnarray*}

That is
\begin{eqnarray}
T(r,f(z))=\overline{N}(r,\frac{1}{f(z)-a})+\overline{N}(r,\frac{1}{f(z)-b})+o(T(r,f)).
\end{eqnarray}

Set
\begin{eqnarray}
\varphi(z)=\frac{L(f(z))(g(z)-f(z))}{(f(z)-a)(f(z)-b)},\\
\psi(z)=\frac{L(g(z))(g(z)-f(z))}{(g(z)-a)(g(z)-b)}.
\end{eqnarray}

If $\varphi(z)\equiv0$, it is a contradiction with $f(z)\not\equiv g(z)$. So $\varphi(z)\not\equiv0$. It is easy to see that $\varphi(z)$ is an entire function. By Lemma 2.1, Lemma 2.5 and Lemma 2.6, we have
\begin{align}
&T(r,\varphi(z))=m(r,\varphi(z))=m(r,\frac{L(f(z))(g(z)-f(z))}{(f(z)-a)(f(z)-b)})\notag\\
\leq & m(r,\frac{L(f(z))f(z)}{(f(z)-a)(f(z)-b)})+m(r,\frac{g(z)}{f(z)}-1)+o(T(r,f))\notag\\
=&o(T(r,f)).
\end{align}

Let $d=a+k(a-b)~(k\neq0,-1)$. Then by Lemma 2.1 and Lemma 2.6 we get
\begin{align}
&m(r,\frac{1}{f(z)-d})=m(r,\frac{L(f(z))(g(z)-f(z))}{\varphi(z)(f(z)-a)(f(z)-b)(f(z)-d)})\notag\\
\leq &m(r,\frac{g(z)}{f(z)}-1)+m(r,\frac{L(f(z))f(z)}{(f(z)-a)(f(z)-b)(f(z)-d)})+o(T(r,f))\notag\\
=&o(T(r,f)).
\end{align}

Set
\begin{align}
F(z)=\frac{f(z)-a}{b-a},\quad G(z)=\frac{g(z)-a}{b-a}.
\end{align}
Because $f(z)$ and $g(z)$ share $a,b$ IM, and $f(z)$ is an entire function of zero-order, we can get $F(z)$ and $G(z)$ are two meromorphic function of zero-order with $\delta(\infty,F(z))=\delta(\infty,G(z))=1$, and $F(z)$ and $G(z)$ share $0,1$ almost IM.

We apply Lemma 2.8 and {\bf Remark 1} to $G(z)$, and by (3.1), we have

\begin{align}
 N(r,\frac{1}{G(z)})+N(r,\frac{1}{G(z)-1})= T(r,f(z))+N(r,\frac{1}{G'(z)})+o(T(r,f)),
\end{align}
which follows from Lemma 2.9 and Lemma 2.10 that
\begin{eqnarray*}
\begin{aligned}
  2T(r,g(z))&=2T(r,G(z))+o(T(r,f))=T(r,f(z))+T(r,G'(z))+\varepsilon_{1}T(r,f)+o(T(r,f))\\
  &=T(r,f(z))+T(r,G(z))+\varepsilon_{1}T(r,f)+o(T(r,f))\\
  &=T(r,f(z))+T(r,g(z))+\varepsilon_{1}T(r,f)+o(T(r,f)),
\end{aligned}
\end{eqnarray*}
which is
\begin{align}
T(r,f(z))= T(r,g(z))+\varepsilon_{1}T(r,f)+o(T(r,f)).
\end{align}

 By the Second Nevanlinna Fundamental Theorem (3.1) and (3.8), we have
\begin{eqnarray*}
\begin{aligned}
&2T(r,f(z))\leq 2T(r,g(z))+\varepsilon_{1}T(r,f)+o(T(r,f))\\
\leq &\overline{N}(r,\frac{1}{g(z)-a})+\overline{N}(r,\frac{1}{g(z)-b})+\overline{N}(r,\frac{1}{g(z)-d})+\varepsilon_{1}T(r,f)+o(T(r,f))\\
\leq &\overline{N}(r,\frac{1}{f(z)-a})+\overline{N}(r,\frac{1}{f(z)-b})+T(r,\frac{1}{g(z)-d})\\
&-m(r,\frac{1}{g(z)-d})+\varepsilon_{1}T(r,f)+o(T(r,f))\\
\leq &T(r,f(z))+T(r,g(z))-m(r,\frac{1}{g(z)-d})+\varepsilon_{1}T(r,f)+o(T(r,f))\\
\leq &2T(r,f(z))-m(r,\frac{1}{g(z)-d})+\varepsilon_{1}T(r,f)+o(T(r,f)).
\end{aligned}
\end{eqnarray*}

Thus
\begin{eqnarray}
m(r,\frac{1}{g(z)-d})=\varepsilon_{1}T(r,f)+o(T(r,f)).
\end{eqnarray}

From the First Fundamental Theorem, Lemma 2.1, Lemma 2.2, (3.5), (3.8), (3.9) and the condition that $f(z)$ is an entire function of zero-order, we obtain
\begin{eqnarray*}
\begin{aligned}
&m(r,\frac{f(z)-d}{g(z)-d})-m(r,\frac{g(z)-d}{f(z)-d})\\
=&T(r,\frac{f(z)-d}{g(z)-d})-N(r,\frac{f(z)-d}{g(z)-d})
-T(r,\frac{g(z)-d}{f(z)-d})+N(r,\frac{g(z)-d}{f(z)-d})\\
=&N(r,\frac{g(z)-d}{f(z)-d})-N(r,\frac{f(z)-d}{g(z)-d})+o(T(r,f))\\
=& N(r,\frac{1}{f(z)-d})-N(r,\frac{1}{g(z)-d})+o(T(r,f))\\
=&T(r,\frac{1}{f(z)-d})-m(r,\frac{1}{f(z)-d})-T(r,\frac{1}{g(z)-d})+m(r,\frac{1}{g(z)-d})+o(T(r,f))\\
=&T(r,f(z))-T(r,g(z))+\varepsilon_{1}T(r,f)+o(T(r,f))=\varepsilon_{1}T(r,f)+o(T(r,f)).
\end{aligned}
\end{eqnarray*}

Thus
\begin{eqnarray}
m(r,\frac{f(z)-d}{g(z)-d})-m(r,\frac{g(z)-d}{f(z)-d})=\varepsilon_{1}T(r,f)+o(T(r,f)).
\end{eqnarray}

It follows from (3.5) and (3.10) that
\begin{align}
&m(r,\frac{f(z)-d}{g(z)-d})=m(r,\frac{g(z)-d}{f(z)-d})+\varepsilon_{1}T(r,f)+o(T(r,f))\notag\\
\leq &m(r,\frac{g(z)-D}{f(z)-d})+m(r,\frac{D-d}{f(z)-d})+\varepsilon_{1}T(r,f)+o(T(r,f))\notag\\
&=\varepsilon_{1}T(r,f)+o(T(r,f)),
\end{align}
where $D=A_{0}(z)d^{(k_{0})}(q_{0}z)+A_{1}(z)d^{(k_{1})}(q_{1}z)+\cdots+A_{1}(z)d^{(k_{j})}(q_{j}z)$.

Rewriting (3.3) we have
$$\psi(z)=[\frac{a-d}{a-b}\frac{L(g(z))}{g(z)-a}-\frac{b-d}{a-b}\frac{L(g(z))}{g(z)-b}][\frac{f(z)-d}{g(z)-d}-1].$$
Then by above and (3.11) we  get
\begin{eqnarray}
T(r,\psi(z))=m(r,\psi(z))+o(T(r,f))=\varepsilon_{1}T(r,f)+o(T(r,f)).
\end{eqnarray}

Now let $m$ and $n$ be two positive integers and let $z_{1}\in S_{(m,n)}(a)\cup S_{(m,n)}(b)$, i.e, $z_{1}$ be a common zero of $f(z)-a$ (resp. $f(z)-b$) and $g(z)-a$ (resp. $g(z)-b$) with multiplicities $m$ and $n$, respectively. (3.2) and (3.3) imply that $n\varphi(z_{1})-m\psi(z_{1})=0$.

Next we consider the following two cases.

{\bf Case1.} \quad $n\varphi(z)- m\psi(z)\equiv0$ for some positive integers $m$ and $n$. It follows that $n\varphi(z)\equiv m\psi(z)$. Then by calculating we have
\begin{eqnarray}
n(\frac{L(f(z))}{f(z)-a}-\frac{L(f(z))}{f(z)-b})\equiv m(\frac{L(g(z))}{g(z)-a}-\frac{L(g(z))}{g(z)-b}),
\end{eqnarray}
which implies that
\begin{eqnarray}
(\frac{f(z)-a}{f(z)-b})^{n}\equiv A(\frac{g(z)-a}{g(z)-b})^{m},
\end{eqnarray}
where $A$ is a nonzero constant.  Hence $n=m$, otherwise we would have a contradiction to (3.8). It follows from (3.14) that

\begin{eqnarray}
B(\frac{f(z)-a}{f(z)-b})\equiv \frac{g(z)-a}{g(z)-b},
\end{eqnarray}
where $B\not= 1$ is a nonzero constant. Thus we have

$$
\frac{b-a}{g(z)-b}=\frac{(B-1)f(z)+(b-aB)}{f(z)-b}.
$$

Since $f(z)$ is an entire function of zero-order, it follows that $f(z)\not= \frac{b-aB}{1-B}.$ Obviously,  $\frac{b-aB}{1-B}\not= a, b.$ Thus we have
\begin{eqnarray*}
\begin{aligned}
2T(r,f(z))&\leq \overline{N}(r,\frac{1}{f(z)-a})+\overline{N}(r,\frac{1}{f(z)-b})+\overline{N}(r,\frac{1}{f(z)-\frac{b-aB}{1-B}})\\
&+o(T(r,f))\le \overline{N}(r,\frac{1}{f(z)-a})+\overline{N}(r,\frac{1}{f(z)-b})+o(T(r,f)),
\end{aligned}
\end{eqnarray*}
which contradicts  (3.1).

\

{\bf Case2.} \quad $n\varphi(z)\not\equiv m\psi(z)$ for any positive integers $m$ and $n$.  Thus we have
\begin{align}
&\overline{N}_{(m,n)}(r,\frac{1}{f(z)-a})+\overline{N}_{(m,n)}(r,\frac{1}{f(z)-b})\leq \overline{N}(r,\frac{1}{n\varphi(z)-m\psi(z)})\notag\\
&\leq T(r,n\varphi(z)-m\psi(z))+o(T(r,f))\notag\\
&\leq T(r,\varphi(z))+T(r,\psi(z))+o(T(r,f))\notag\\
&=\varepsilon_{1}T(r,f)+o(T(r,f)),
\end{align}
for all positive integers $m$ and $n$.

Thus by (3.8) and (3.16), we  get
\begin{align}
&T(r,f(z))\leq \overline{N}(r,\frac{1}{f(z)-a})+\overline{N}(r,\frac{1}{f(z)-b})+o(T(r,f))&\notag\\
\leq &\overline{N}_{1)}(r,\frac{1}{f(z)-a})+\overline{N}_{2}(r,\frac{1}{f(z)-a})+\overline{N}_{3}(r,\frac{1}{f(z)-a})+\overline{N}_{4}(r,\frac{1}{f(z)-a})\notag\\
+ &\overline{N}_{(5}(r,\frac{1}{f(z)-a})+\overline{N}_{1)}(r,\frac{1}{f(z)-b})+\overline{N}_{2}(r,\frac{1}{f(z)-b})+\overline{N}_{3}(r,\frac{1}{f(z)-b})\notag\\
+ &\overline{N}_{4}(r,\frac{1}{f(z)-b})+\overline{N}_{(5}(r,\frac{1}{f(z)-b})+o(T(r,f))\notag\\
\leq &\sum_{n=1}^{4}\sum_{m=1}^{4}\overline{N}_{(m,n)}(r,\frac{1}{f(z)-a})+\overline{N}_{(5}(r,\frac{1}{g(z)-a})+\overline{N}_{(5}(r,\frac{1}{f(z)-a})\notag\\
+&\sum_{n=1}^{4}\sum_{m=1}^{4}\overline{N}_{(m,n)}(r,\frac{1}{f(z)-b})+\overline{N}_{(5}(r,\frac{1}{g(z)-b})+\overline{N}_{(5}(r,\frac{1}{f(z)-b})+o(T(r,f))\notag\\
\leq &\frac{1}{5}[N(r,\frac{1}{f(z)-a})+N(r,\frac{1}{f(z)-b})]+\frac{1}{5}[N(r,\frac{1}{g(z)-a})+N(r,\frac{1}{g(z)-b})]+16\varepsilon_{1}T(r,f)+o(T(r,f))\notag\\
\leq &\frac{2}{5}T(r,f(z))+\frac{2}{5}T(r,g(z))+\varepsilon_{1}T(r,f)+16\varepsilon_{1}T(r,f)+o(T(r,f))\notag\\
= &\frac{4}{5}T(r,f(z))+16\varepsilon_{1}T(r,f)+o(T(r,f)),
\end{align}
it follows from above that
\begin{align}
(\frac{1}{5}-16\varepsilon_{1})T(r,f(z))=o(T(r,f)).
\end{align}
We take $\varepsilon<\frac{1}{80}$ in Lemma 2.8, and thus we obtain from (3.18) that $T(r,f(z))=o(T(r,f))$, a contradiction.

This completes the proof of Theorem 1.

\

{\bf Acknowledgements} The author would like to thank to anonymous referees for their helpful comments.


\end{document}